%Paper: dg-ga/9503011
%From: Vladimir.Pestov@vuw.ac.nz (V. Pestov)
%Date: Tue, 28 Mar 1995 11:44:04 +1300

\documentstyle{amsppt}
\magnification=\magstep 1
\TagsOnRight
\NoBlackBoxes
\leftheadtext{V. Pestov}
% macros
\def\norm #1{{\left\Vert\,#1\,\right\Vert}}

\def\R {{\Bbb R}}
\def\C {{\Bbb C}}
\def\N{{\Bbb N}}

\def\Lip{{\operatorname{Lip}\,}}

\def\QED{\nobreak\quad\ifmmode\roman{Q.E.D.}\else{\rm Q.E.D.}\fi}
\def\ml #1\endml{\email #1\endemail}
\long\def\block #1\endblock{\vskip 6pt
        {\leftskip=3pc \rightskip=\leftskip
        \noindent #1\endgraf}\vskip 6pt}
\long\def\ext #1\endext{\block #1\endblock}

\def\tri
{\ifhmode\unskip\nobreak\fi\ifmmode\ifinner\else\hskip 5pt\fi\fi
 \hfill\hbox{\hskip 5pt$\blacktriangle$\hskip 1pt}}
% end of macros
% Counters for references
\newcount\refAdams
\newcount\refArh
\newcount\refAP
\newcount\refBlum
\newcount\refBiS
\newcount\refBiSi
\newcount\refD
\newcount\refDo
\newcount\refEe
\newcount\refE
\newcount\refG
\newcount\refKG
\newcount\refP
\newcount\refRam
\newcount\refRa
\newcount\refRu

\newcount\refno
\refno=0

\advance\refno by 1\refAdams=\refno
\advance\refno by 1\refArh=\refno
\advance\refno by 1\refAP=\refno
\advance\refno by 1\refBlum=\refno
\advance\refno by 1\refBiS=\refno
\advance\refno by 1\refBiSi=\refno
\advance\refno by 1\refD=\refno
\advance\refno by 1\refDo=\refno
\advance\refno by 1\refEe=\refno
\advance\refno by 1\refE=\refno
\advance\refno by 1\refG=\refno
\advance\refno by 1\refKG=\refno
\advance\refno by 1\refP=\refno
\advance\refno by 1\refRam=\refno
\advance\refno by 1\refRa=\refno
\advance\refno by 1\refRu=\refno

%end of counters for reference
\topmatter
\title
Analytic subsets of Hilbert spaces
$\dag$
\endtitle
\author
Vladimir Pestov
\endauthor
\address
Department of Mathematics,
Victoria University of Wellington,  P.O. Box 600, Wellington,  New
Zealand
\endaddress
\email {\tt vladimir.pestov$\@$vuw.ac.nz}
\endemail
\abstract{We show that every complete metric space is homeomorphic to
the locus of zeros of an entire
analytic map from a complex Hilbert space to a complex Banach space.
As a corollary, every separable complete metric space is homeomorphic
to the locus of zeros of an entire analytic map between two complex
Hilbert spaces.}

\endabstract
\endsubjclass
\endtopmatter
\document
\footnote""{$\dag$ Research Report RP-95-157,
Department of Mathematics,
Victoria University of Wellington, March 1995.}

\subheading{\S 1}
Douady had observed  \cite{\the\refDo}
that every compact metric space is
homeomorphic to the locus of zeros of an analytic (in fact,
continuous polynomial) map between two suitable Banach spaces.
In other words, every compact metric space is homeomorphic to an
algebraic subset of a complex Banach space.
Proving Douady's conjecture, the present author has shown
\cite{\the\refP} that every complete metric space is
isometric to
an algebraic subset of a complex Banach space.
Hilbert spaces provide a particularly
favourable setting for Banach analytic geometry
(cf. \cite{\the\refRu}), and
the following question suggested by Norm
Dancer is very natural: what can be said about the topology of
analytic subsets of Hilbert (or just reflexive Banach) spaces?

A previously known result belongs to
Ramis \cite{\the\refRa}
who embedded the Cantor set topologically in a Hilbert space
as an algebraic subset.

\proclaim{Main Theorem}
Every
complete metric space is homeomorphic to the
precise locus of zeros of
an entire analytic map from a complex Hilbert space to a
complex Banach space.
\endproclaim

Call a Banach analytic space \cite{\the\refD}
an {\it Hilbert analytic space}
if it is modelled on analytic subsets of an Hilbert space
(zeros of Banach-valued analytic maps).
The following corollary betters both Douady's
and the present author's earlier results.

\proclaim{Corollary 1} A paracompact topological space admits the
structure of an Hilbert analytic space if and only if it is metrizable
with a complete metric.
\endproclaim

\proclaim{Corollary 2}
Every separable
complete metric space is homeomorphic to the
locus of zeros of
an entire map between two separable complex Hilbert spaces.
\endproclaim

Every closed subset of a separable real Hilbert space
is the precise locus of zeros of a $C^\infty$ functional
(see e.g. \cite{\the\refEe}, 2.C), but
obviously this result does not extend to real
analytic functionals. However, one has:

\proclaim{Corollary 3}
Every separable
complete metric space is homeomorphic to the locus of zeros of
an
entire functional on a separable real Hilbert space.
\endproclaim

The latter result indicates that theory of
finitely defined analytic subsets of Banach (even Hilbert) manifolds
\cite{\the\refRam}, \cite{\the\refRa}
only has topological substance in a complex case.

\definition{Remark 1}
One can hardly expect the above Main Theorem
to admit an isometric version.
Already the
unit segment $I=[0,1]$ with the usual Euclidean
distance cannot be isometrically embedded as
an {\it analytic} subset in a Hilbert space. If
$I\hookrightarrow\Cal H$ is such an isometry, then $I$ forms a segment of a
{\it metric line} in $\Cal H$ and therefore a segment
of a real affine line (cf. \S 3 in \cite{\the\refBlum}).
As an easy consequence, the locus of zeros of an appropriate
non-constant real analytic map on an interval in
$\R$ contains an open non-void
subinterval of $I$, which is impossible.
\enddefinition

\subheading{\S 2}
All Banach spaces and algebras below are complex unless otherwise
stated. In particular, the Sobolev spaces $W^{n,m}(\Omega)$
\cite{\the\refAdams},
the spaces $C^n[a,b]$,
and the Banach algebra $\Lip [0,1]$ of Lipschitz functions
are formed by complex-valued functions.
We denote the locus of zeros of a mapping $f$ by $V(f)$.
An $l_p$-sum of a family of Banach spaces is
denoted by $\oplus^{l_p}$. The {\it diagonal product}
$\Delta_\alpha f_\alpha$ of a family of maps
$f_\alpha\colon X\to Y_\alpha$ sends an $x\in X$ to
$(f_\alpha(x))\in\prod_\alpha X_\alpha$.

\proclaim{Theorem 1 {\rm (Douady \cite{\the\refDo};
\cite{\the\refRa}, 1.3.A)}}
The maximal ideal space, $X$, of a unital Banach algebra $A$ is
the locus of zeros of a continuous $2$-polynomial on
the Banach dual $A^\prime$.
\qed\endproclaim

A Sobolev space $W^{2,2}(0,1)$ forms a separable
Hilbert space and also a commutative unital
Banach algebra
under the pointwise multiplication
(\cite{\the\refAdams}, 3.5 and 5.23).

\proclaim{Assertion 1} The maximal ideal space of $W^{2,2}(0,1)$
is canonically homeomorphic to the unit segment $[0,1]$ in both
the weak topology and the strong dual topology.
\endproclaim

\demo{Proof} The canonical embedding $C^2[0,1]\to C^1[0,1]$ is
the composition of
an obvious continuous embedding
$i_2\colon C^2[0,1]\to W^{2,2}(0,1)$ and a continuous embedding
$i_1\colon W^{2,2}(0,1)\to C^1[0,1]$ asserted by the Sobolev Embedding
Theorem (\cite{\the\refAdams}, 5.4.I.C), as can be seen from
\cite{\the\refAdams}, 5.2.
 Therefore, both $i_1$ and $i_2$ are
homomorphisms of unital Banach algebras.
Since the elements of
maximal ideal spaces of both $C^1[0,1]$ and $C^2[0,1]$
are in a natural one-to-one correspondence with the
points of $[0,1]$
(\cite{\the\refG}, 1.3.3),
the same is true for $W^{2,2}(0,1)$.
The Gelfand space of
$\Lip [0,1]$
is canonically homeomorphic
to $[0,1]$ in both the strong and the weak topologies
\cite{\the\refRa}, \cite{\the\refDo},
\cite{\the\refP}.
It follows that the dual operator to the composition of
$i_1$ and the canonical
homomorphism $C^1[0,1]\to\Lip [0,1]$
maps the maximal ideal space of $\Lip [0,1]$
onto that of $W^{2,2}(0,1)$ in a one-to-one fashion and continuously
with respect to both pairs of topologies.
\qed\enddemo

\proclaim{Assertion 2} A separable Hilbert space $\Cal H$ contains
a topological copy of the unit interval,
$I$, as the locus of zeros of a continuous
$2$-polynomial on $\Cal H$ in such a way that
the left endpoint $0_I\in I$ is
the zero element of $\Cal H$, and every closed subset of $I$
is the intersection of $I$ with the locus of zeros of a continuous
$1$-polynomial.
\endproclaim

\demo{Proof} Denote by $\Cal H$ the dual space of distributions
$W^{2,2}(0,1)^\prime$.
In view of the above two results,
the translation, $X-0_I$, of the maximal ideal
space, $X\cong I$, of $W^{2,2}(0,1)$
is the precise locus of
zeros of a continuous $2$-polynomial on $\Cal H$.
Every
closed subset, $F$, of $I\cong X$ is the locus of zeros of a
suitable $C^\infty$-function, $f$, on $I$
(e.g. \cite{\the\refEe}, 2.C). Viewed as an element of
$C^2[0,1]\subset W^{2,2}(0,1)\cong\Cal H^\prime$, this function determines a
continuous linear functional $\bar f$ on $\Cal H$ with
the property $F=X\cap V(\bar f)$. Finally,
$F-0_I=(X-0_I)\cap V(\bar f-\bar f(0_I))$.
\qed\enddemo

\proclaim{\S 3. Lemma 1}
The intersection of a family of subsets of a
Banach space $E$, each of which is the locus of zeros of a continuous
(resp. continuous homogeneous)
polynomial map of degree $\leq n$ (resp. $n$)
is the locus of zeros of  a continuous (resp. continuous homogeneous)
polynomial map of degree $\leq n$
(resp. $n$).
\endproclaim

\demo{Proof}
Let $p_\alpha\colon E\to F_\alpha$ be continuous
Banach space valued $n$-polynomial maps.
One can assume without loss in generality that the norm
\cite{\the\refBiS}, \cite{\the\refD}
of every
homogeneous component of each polynomial $p_\alpha$ is $\leq 1$.
The map
$p\equiv\Delta_\alpha p_\alpha\colon E\to \oplus^{l_\infty} F_\alpha$
is a well-defined continuous polynomial of
degree $\leq n$. Indeed, $p$ can be represented as the sum of
diagonal products $p^{(i)}\equiv
\Delta_\alpha p_\alpha^{(i)}$ of homogeneous components of $p_\alpha=
\sum_{i=0}^np_\alpha^{(i)}$. Since for each $t\in\C$ one has
$p^{(i)}(tx)=t^ip^{(i)}(x)$, each $p^{(i)}$ is a homogeneous polynomial
of degree $i$ (\cite{\the\refBiS}, coroll. 3.H).
Furthermore, the norm of every homogeneous
component of $p$ does not exceed $1$, and therefore
every $p^{(i)}$ is continuous
(\cite{\the\refBiS}, Th. 1; \cite{\the\refD}, prop. 1),
as well as $p$ itself. It is clear that
$\cap_\alpha V(p_\alpha)=V(p)$. Finally, if each $p_\alpha$ is
homogeneous of degree $n$, then so is $p$.
\qed\enddemo

\proclaim{Lemma 2}
Let $\Cal H_\alpha$ be a family of Hilbert spaces.
Then the union $\cup_\alpha \Cal H_\alpha$ of these spaces canonically
embedded into their Hilbert direct sum
$\Cal H=\oplus^{l_2} \Cal H_\alpha$
forms the locus of zeros of a continuous
homogeneous polynomial map on $\Cal H$ of degree $2$.
\endproclaim

\demo{Proof} For each pair of indices $\alpha,\beta$,
$\alpha\neq\beta$, let $\pi_{\alpha,\beta}$ be the projection of
$H$ to $H_\alpha\oplus H_\beta$, and let $i_{\alpha,\beta}$ stand
for the canonical embedding of $H_\alpha\times H_\beta$ into
$H_\alpha\hat\otimes H_\beta$ (the projective
tensor product). The locus of zeros of the continuous homogeneous
$2$-polynomial $i_{\alpha,\beta}\circ \pi_{\alpha,\beta}$
consists of such elements $(x_\gamma)\in H$ that
at most one of the two
coordinates $x_\alpha,x_\beta$ is non-vanishing.
The intersection of sets
$V(i_{\alpha,\beta}\circ \pi_{\alpha,\beta})$ over all pairs
$(\alpha,\beta)$ consists of all $(x_\gamma)\in H$ with at most one
coordinate non-vanishing, as desired. Now apply Lemma 1.
\qed\enddemo

\proclaim{Lemma 3}
The intersection of a family of subsets of a
Banach space $E$, each of which is the locus of zeros of a continuous
polynomial map, is the locus of zeros of an entire map on $E$.
\endproclaim

\demo{Proof}
In view of Lemma 1, grouping together polynomials of the same
degree, one can assume
the family of polynomials to be
countable, with $p_n\colon E\to F_n$ and
$\operatorname{deg}\,p_1<\operatorname{deg}\,p_2<\dots$.
One can also assume, multiplying $p_n$ by a suitable
scalar if necessary, that for every $n$ and every
$i=0,1,\dots, \operatorname{deg}\,p_n$ one has
$\sup \{\norm{p_n^{(i)}(x)}\colon \norm x\leq i\}\leq 2^{-i}$.
The diagonal product
$p^{(i)}=\Delta_np_n^{(i)}\colon E\to \oplus^{l_\infty}_nF_n$
is a homogeneous polynomial map
with
$\sup \{\norm{p^{(i)}(x)}\colon \norm x\leq i\}\leq 2^{-i}$,
and the series
$\sum_{i=0}^\infty p^{(i)}$ converges uniformly on every bounded set
in $E$ to the mapping
$p=\Delta p_n\colon E\to \oplus^{l_\infty}_nF_n$.
Therefore, $p$ is an entire map
(\cite{\the\refBiSi}, Sect. 8), and also
clearly $V(p)=\cap_nV(p_n)$.
\qed\enddemo

\definition{Remark 2}
An analysis of the Example
in \cite{\the\refRa}, p. 84 shows that
the intersection of an infinite
family of algebraic sets is not in general an algebraic set.
\enddefinition

\proclaim{Lemma 4 {\rm (cf. \cite{\the\refRa}, prop. II.1.1.1)}}
The union of finitely many subsets of a Banach space, each of which is
the locus of zeros of a continuous
polynomial map, is the locus of zeros of a continuous
polynomial map.
\endproclaim

\demo{Proof}
If $p_i,~i=1,\dots, n$ are continuous polynomial maps $E\to F_i$,
then $\cap_{i=1}^n V(f_i)=V(f)$, where $f$ is the
composition of the diagonal product
$\Delta f_i\colon x\mapsto (f_1(x), f_2(x),\dots,$
$f_n(x))\in F_1\times F_2\times\dots\times F_n$
and the canonical $n$-linear mapping
$F_1\times F_2\times\dots\times F_n\to
F_1\hat\otimes F_2\hat\otimes\dots\hat\otimes F_n$.
\qed\enddemo

\proclaim{Lemma 5}
Let $\Cal H_\alpha$ be a family of
Hilbert spaces, let $k\in\N$, and let, for each $\alpha$,
$X_\alpha$ be the locus of zeros of
a continuous Banach space valued $k$-polynomial on $\Cal H_\alpha$
such that $0\in X_\alpha$. The
union
$\cup_\alpha X_\alpha$ of the sets $X_\alpha$ canonically embedded
into the $l_2$-sum $\Cal H=\oplus^{l_2}_\alpha \Cal H_\alpha$ is
the locus of zeros of
a continuous Banach space valued $\max\{k,2\}$-polynomial.
\endproclaim

\demo{Proof} Let $X_\alpha=V(p_\alpha)$, where
$p_\alpha\colon \Cal H_\alpha\to E_\alpha$ are continuous
$k$-polynomials. Then $\cup_\alpha X_\alpha$ is the intersection of
$\cup_\alpha \Cal H_\alpha$ and
the intersection of the family of
loci of zeros of homogenuous $k$-polynomials
of the form $p_\alpha\circ\pi_\alpha$, where
$\pi_\alpha\colon \Cal H\to \Cal H_\alpha$ is the projection.
Now it remains to apply Lemmas 2 and 1.
\qed\enddemo

The symbol $J(\frak m)$  denotes the ``metrizable hedgehog
of spininess $\frak m$,'' that is, a bouquet of $\frak m$
copies of the unit interval $I$ with a common basepoint $0_I$,
endowed with the maximal metric inducing the usual distance on
each copy of $I$.
(\cite{\the\refE}, Ch. 4; \cite{\the\refAP}, ch. 2, problem 251).

\proclaim{Lemma 6} Let $\frak m$ be a cardinal.
The Hilbert space $\Cal H$ of weight $\frak m$ contains
a bounded subset homeomorphic to $J(\frak m)$
in such a way that $J(\frak m)$ itself, and
every element of a certain
basis of closed subsets of this space, are the loci of zeros
of suitable Banach space-valued
continuous $2$-polynomials on $\Cal H$.
\endproclaim

\demo{Proof} Let $I_\alpha$, $\alpha<\frak m$, be copies of the unit
interval embedded in a separable Hilbert space $\Cal H_\alpha$
as in Assertion 2. According to Lemma 5, the union
$\cup_{\alpha<\frak m}I_\alpha$ formed in a Hilbert space
$\Cal H=\oplus^{l_2}_{\alpha<\frak m}\Cal H_\alpha$
is the locus of zeros
of a Banach space-valued continuous $2$-polynomial.
It is easy to verify that
$\cup_{\alpha<\frak m}I_\alpha$ is bounded and
homeomorphic to the metric hedgehog.
If a closed subset $F\subseteq J(\frak m)$ contains zero,
then it can be represented as the union of closed subsets
$F_\alpha\subseteq I_\alpha$ each of which contains zero, and
by force of Lemma 5 $F$ is
the locus of zeros
of a continuous $2$-polynomial map.
It suffices now to check the required property for a closed subset
of the form $\cup_\alpha F_\alpha$ where each $F_\alpha$ is
a copy of the same
closed subinterval $[a,1]\subset I_\alpha$,
$a>0$. Let $f\colon \Cal H\to\C$ be a linear
bounded functional with $f^{-1}(\{1\})=[a,1]\cap I$ (Assertion 2).
The operator
$f^\frak m\colon \oplus^{l_2}\Cal H\to l_\infty(\frak m)$
of the form
$(x_\alpha)\mapsto (f(x_\alpha))$ is bounded
linear,  and
$F$ is the intersection of $J(\frak m)$ with the inverse image
of the constant sequence $(1)\in l_\infty(\frak m)$
under $f^\frak m$.
\qed\enddemo

\proclaim{Lemma 7} Let $\frak m$ be a cardinal.
The Hilbert space $\Cal H$ of weight $\frak m$ contains
a subset homeomorphic to the countable
Tychonoff power $J(\frak m)^{\aleph_0}$
in such a way that $J(\frak m)^{\aleph_0}$ itself, and
every element of a certain
basis of closed subsets of this space, are the loci of zeros
of suitable Banach space-valued
continuous polynomials on $\Cal H$.
\endproclaim

\demo{Proof} Let $i\colon J(\frak m)\to \Cal H$
be a homeomorphic embedding from Lemma 6.
A map from
$J(\frak m)^{\aleph_0}$ into the
$l_2$-sum of countably many copies of
$\Cal H$, given by the rule
$(x_n)\mapsto (2^{-n}i(x_n))$, is a homeomorphic
embedding by force of
a simple argument based upon the boundedness of $i(J(\frak m))$
and identical to the familiar proof of the fact that the Hilbert cube,
$Q^\infty$ (see e.g. \cite{\the\refAP}, 3.4),
is homeomorphic to $I^{\aleph_0}$.
Applying Lemma 1 to the compositions of
the projections
$\oplus^{l_2}_n\Cal H_n\to\Cal H_n$ with
the $2$-polynomial maps
determining subsets $2^{-n}i(J(\frak m))$ in the factors
$\Cal H_n\cong\Cal H$, one concludes that
$J(\frak m)^{\aleph_0}$
the locus of zeros of a continuous $2$-polynomial Banach space
valued map. In view of Lemma 6, the complement to each
element of a standard open subbase for
$J(\frak m)^{\aleph_0}$ (a cylinder over a basic open subset of
$J(\frak m)$) is the locus of zeros of a continuous
$2$-polynomial. Finally,
the complement to each standard basic
open subset in
$J(\frak m)^{\aleph_0}$, being a finite union of complements to
standard subbasic elements, forms by force of Lemma 4
the
locus of zeros of a continuous polynomial map.
\qed\enddemo

\subheading{\S 4. Proof of the Main Theorem}
It follows from Lemma 7 and Lemma 3
that every closed subset of
$J(\frak m)^{\aleph_0}$ is homeomorphic to the locus of zeros
of a suitable Banach space valued entire map on a Hilbert
space $\Cal H$ of weight $\frak m$: indeed, every such subset is an
intersection of a family of loci of zeros of continuous
polynomial maps on $\Cal H$.
But every completely metrizable space of weight
$\leq\frak m$ is homeomorphic to a closed
subspace of $J=J(\frak m)^{\aleph_0}$
(\cite{\the\refE}, Problem 4.4.B.)\qed

\subheading{Proof of Corollary 1}
A paracompact space, locally metrizable with a
complete metric, is itself metrizable with a complete metric
\cite{\the\refArh}.\qed

\subheading{Proof of Corollary 2} The Main Theorem enables one
to represent any separable
completely metrizable space as the locus of zeros of
an entire map from a separable
complex Hilbert space to a separable
complex Banach space, $E$. It remains to
compose this map with an isometric embedding of
$E$ into $C[0,1]$ (see e.g. \cite{\the\refKG}, Th. 24)
and a canonical linear contraction
$C[0,1]\to\Cal L^2(0,1)$.\qed

\subheading{Proof of Corollary 3}
A map from Corollary 2,
viewed as an entire real analytic map between two real
Hilbert spaces, should be composed
with the norm on the image
Hilbert space, which in the real case is a
continuous $2$-polynomial functional. \qed

\subheading{\S 5. Questions. 1}
Is every complete metric space
homeomorphic to the locus of zeros of an analytic map
between two Hilbert spaces?
\smallskip
{\bf 2.}
Is every complete metric space
homeomorphic to an {\it algebraic} subset of a Hilbert space?
\smallskip
{\bf 3.} Let $E$ be an infinite-dimensional Banach space of
weight $\frak m$. Is every complete metric space of weight
$\leq\frak m$ homeomorphic to an analytic subset of $E$?

\subheading{\S 6. Acknowledgment} The author is grateful to
Prof. E.N. Dancer
(The University of Sydney) whose question, asked
after author's talk
at the Armidale OzMS-94 Annual Conference, has stimulated the present
investigation.

\enddocument

\Refs
\widestnumber\key{Fl2}
\vskip0.3truecm

\ref\key\the\refAdams
\by R\.A\. Adams
\book Sobolev Spaces
\publ Academic Press
\publaddr NY--San Fransisco--London
\yr 1975
\endref

\ref\key \the\refArh
\by
A\.V\. Arhangel'ski\u\i
\paper \u Cech-complete topological spaces
\jour Vestnik Moskov\. Univ\. Ser\. Mat\.
\vol 2
\yr 1961
\pages 37--40
\lang Russian
\endref

\ref\key\the\refAP
\by A\.V\. Arhangel'ski\u\i\ and V\.I\. Ponomarev
\book Fundamentals of General Topology: Problems and Exercises
\publ D\. Reidel Pub\. Co\.
\publaddr Dordrecht, Boston
\yr 1984
\endref

\ref\key\the\refBlum
\by L\.M\. Blumenthal
\paper Four-point properties and norm postulates
\inbook The Geometry of Metric and Linear Spaces
(Proceedings, Michigan 1974), Lecture Notes in Mathematics
{\bf 490}
\publ Springer-Verlag
\publaddr Berlin--Heidelberg--NY
\pages 1--13
\endref

\ref\key\the\refBiS
\by J\. Bochnak and J\. Siciak
\paper Polynomials and multilinear mappings in topological
vector spaces
\jour Studia Mathematica
\vol 39
\yr 1971
\pages 59--76
\endref

\ref\key\the\refBiSi
\bysame
\paper Analytic functions in topological vector spaces
\jour ibidem
\vol 39
\yr 1971
\pages 77--112
\endref

\ref\key \the\refD
\by A\. Douady
\paper Le probl\`eme des modules pour les sous-espaces analytiques
compacts d'un espace analytique donn\'e
\jour Ann\. Inst\. Fourier, Grenoble
\yr 1966
\vol 16
\pages 1--95
\endref

\ref\key \the\refDo
\bysame
\paper A remark on Banach analytic spaces
\inbook Symposium on Infinite-Dimensional Topology
(Annals of Mathematics Studies, {\bf 69})
\publ Princeton University Press and University of Tokyo Press
\publaddr Princeton, NJ
\yr 1972
\pages 41--42
\endref

\ref\key \the\refEe
\by J\. Eells,  jr\.
\paper A setting for global analysis
\jour Bull\. Amer\. Math\. Soc\.
\vol 72
\yr 1966
\pages 751--807
\endref

\ref\key\the\refE
\by R. Engelking
\book General Topology
\publ PWN -- Polish Scientific Publishers
\publaddr Warszawa
\yr 1977
\endref

\ref\key
\the\refG
\by T\.W\. Gamelin
\book Uniform Algebras
\publ Prentice-Hall
\publaddr Englewood Cliffs, N.J.
\yr 1969
\endref

\ref\key\the\refKG
\by A\.A\. Kirillov and A\.D\. Gvishiani
\book Theorems and Problems in Functional Analysis
\publ Springer-Verlag
\bookinfo translated from Russian by Harold H\. McFaden
\yr 1982
\publaddr NY
\endref

\ref\key \the\refP
\by V\. Pestov
\paper Douady's conjecture on Banach analytic spaces
\jour C\.R\. Acad\. Sci\. Paris, S\'er\. I
\yr 1994
\vol 319
\pages 1043--1048
\endref

\ref\key\the\refRam
\by J\.-P\. Ramis
\paper Sous-ensembles analytiques d'une vari\'et\'e
analytique banachique complexe
\jour C\.R\. Acad\. Sci\. Paris, S\'er\. A
\yr 1968
\vol 267
\pages 732--735
\endref

\ref\key\the\refRa
\by J\.-P\. Ramis
\book Sous-ensembles analytiques d'une vari\'et\'e banachique
complexe
\publ Springer-Verlag
\publaddr Berlin--Heidelberg--NY
\yr 1970
\endref

\ref\key\the\refRu
\by G. Ruget
\paper
A propos des cycles analytiques de dimension infinie
\jour Inventiones Math.
\vol 8
\yr 1969
\pages 267--312
\endref


\endRefs
\bye